\documentclass[12pt]{article}      
\usepackage{amssymb}
\newcommand{\got}{\mathfrak}
\begin{document} \sloppy
\title{Addendum to: A simultaneous Frobenius splitting for closures of conjugacy classes of
nilpotent matrices, by V. B. Mehta and Wilberd van der Kallen}  
\author{Wilberd van der Kallen}   
\maketitle       

\section{Introduction}
Theorem 3.8 of \cite{MehtavdK} implies that a certain Frobenius splitting of 
$G\times^B\got b$ is compatible with $G\times^B\got n_P$, where $\got n_P$
denotes the Lie algebra of the unipotent radical of a certain
kind of standard parabolic subgroup $P$ of $G=Gl_n$.
In Exercise 5.1.E.6 of their book \cite{BrionKumar} Brion and Kumar ask to
prove this same fact for any standard parabolic subgroup.
Their comments 5.C suggest that they thought this has already been done in 
\cite{MehtavdK}. However, in \cite{MehtavdK} we did not need all standard
parabolic subgroups and we only treated a class that is slightly easier.
Let us now do the exercise, by discussing the necessary modifications.

\subsection{A partial order.}\label{poset}
We simplify the partial order of \cite[3.1]{MehtavdK}.
We put a partial
order on the set $
{\cal I}= [1,n]\times [1,n]$ which indexes the coordinates on
$\got g$. We declare that
$$(i,j)\leq (r,s) \Longleftrightarrow (i\geq r \mbox{ and } j\leq s)
$$
If $S$ is an ideal for this partial order, i.e.\ if $(i,j)\leq (r,s)$
and $(r,s)\in S$ imply $(i,j)\in S$, then we define ${\got b}[S]$ to be
the subspace of $\got b$ consisting of the matrices $X$ with $X_{ij}=0$
for $(i,j)\in S$. One easily sees that such a subspace is an
ideal, and all ideals of $\got b$ arise this way. Note that
ideals in $\got b$ are $B$ invariant.
Let us agree to use the notation
${\got b}[S]$ only when $S$ is an ideal for the partial order.
We will find a
Frobenius splitting for all $G\times ^B {\got b}[S]$ simultaneously.

We must argue a little differently than in \cite[3.7]{MehtavdK}.
In particular, we can not use \cite[Lemma 3.3]{MehtavdK} now.
\subsection{Start of proof.}
We argue by induction on the size of $S$ to show that specialization leads
to the formulas indicated in \cite[3.7]{MehtavdK},
but we will go in the other direction to
prove that one has Frobenius splittings. The formula for $\sigma[S]$ is by
definition correct when $\got b[S]$ equals $\got b$. (Note that in this
case $i>j$ for $(i,j)\in S$ so that $\delta_r[S]$ vanishes for $r\leq n$.)
Therefore let us now assume $S$ contains a maximal element $(s,t)$ with
$s\leq t$. We assume the formulas true for $S'=S-\{(s,t)\}$. Put $r=s+n-t$.
For $(g,X)\in U^-\times \got b[S']$
we first claim that with $M=X+\delta_r[S']$, the determinant
$\det{((gMg^{-1})_{\leq r,\leq r})}$, which is of degree one in $X_{st}$, 
is divisible by $X_{st}$.
Now this can be checked by putting $X_{st}$ equal to zero and showing that the 
rank of $M_{\leq r,\leq n}$ becomes strictly less than $r$.
Indeed 
$M_{\leq r, \leq n}$ is a block matrix
$$\pmatrix{\alpha&\beta&\gamma\cr 0&X_{st}&\delta\cr 0& 0&\epsilon}$$
which becomes $$\pmatrix{\alpha&\beta&\gamma\cr 0&0&\delta\cr 0& 0&\epsilon}$$
when you put $X_{st}$ equal to zero.
The submatrix $$\pmatrix{\alpha&\beta}$$
has rank at most $s-1$ and the submatrix $$\pmatrix{\gamma\cr \delta\cr \epsilon}$$
has rank at most $n-t$, so together the rank is at most $r-1$ indeed.
 We may use $X_{st}$ as
the $f$ of \cite[3.5]{MehtavdK}, at least over the open subset $U^-$ of
$G/B$. As $U^-\times \got b[S]$ is dense in $G\times^B \got b[S]$, the
hypotheses for the residue construction are satisfied and we only need to
check that it replaces the factor $\det((g(X+\delta_r[S'])g^{-1})_{\leq r,\leq r})$ in
the product for $\sigma[S']$ by the factor
$\det((g(X+\delta_r[S])g^{-1})_{\leq r,\leq r})$. Indeed
one must put $X_{st}$ equal to zero in the regular function
$$\det((g(X+\delta_r[S'])g^{-1})_{\leq r,\leq r})/X_{st}$$
And this gives the same as putting $X_{st}$ equal
to zero in $\det((g(X+\delta_r[S])g^{-1})_{\leq r,\leq r}$.
The rest of the proof proceeds as before.

Note that we are dealing here with a residually normal crossing situation in the sense
of \cite{MPL}.


\begin{thebibliography}{}
\bibitem{BrionKumar}M. Brion and S. Kumar,
Frobenius Splitting Methods in Geometry and Representation Theory, 
Birkh\"auser Boston 2005.
\bibitem{MehtavdK}V.B. Mehta and W. van der Kallen,
A simultaneous Frobenius splitting for closures of conjugacy classes of
nilpotent matrices, Compositio Math.\ 84 (1992), 211--221.
\bibitem{MPL}V. Lakshmibai, V.B. Mehta and A.J. Parameswaran,
Frobenius splittings and blowups, J. Algebra, 208 (1998), 101--128. 
\end{thebibliography}
\end{document}